\documentclass[11pt]{amsart}
\usepackage[english]{babel}
\usepackage{amssymb}
\usepackage{amsmath}
\newcommand{\ignore}[1]{}

\newtheorem{theorem}{Theorem}

\newtheorem{corollary}{Corollary}

\theoremstyle{definition}

\newtheorem{remark}{Remark}
\newtheorem{example}{Example}



\title{A condition for cyclic algebras to be split}
\author{S. Pumpl\"un}

\email{ susanne.pumpluen@nottingham.ac.uk}
\address{School of Mathematics\\
University of Nottingham\\
University Park\\
Nottingham NG7 2RD\\
United Kingdom
}

\begin{document}

\maketitle

\begin{abstract} A simple sufficient condition for certain cyclic algebras of odd degree $d$ to be split is presented.
It employs certain binary forms of degree $d$ and the values they represent. A similar sufficient condition for
certain Albert algebras not to be division algebras
is found as well.
\end{abstract}

\section{Introduction}

For a quaternion algebra $(a,b)_k$ over a field $k$ of characteristic not two,
$(a,b)_k\cong\mbox{Mat}_2(k)$ if and  only if the binary quadratic form $\langle a,b \rangle$
represents $1$ [L, Theorem 2.7(7), p.~58].
Given a cyclic central simple (associative) algebra  $A=(l,a)$ over $k$ of degree $n$, it is well known that
  $$ A\cong\mbox{Mat}_n(k) \text{ if and only if } a\in n_{l/k}(l^\times)$$
   (cf. for instance [KMRT, (30.6), p.~415]).
This criterion translates the question whether $A$ is split into the problem whether a certain element
is represented by a homogeneous form (the reduced norm of $l/k$) of degree $d$ in $d$ indeterminates.
In order to give a sufficient criterion for $A$ to be split it is sometimes enough, however, to check a much simpler
- binary - form of degree $d$ which arises by restricting the reduced norm of the algebra $A$ to a suitable subspace of $A$:
 let $k$ be a field of characteristic not 3 which contains a primitive
third root of unity, and let $A\cong (l,a)$ be a cyclic algebra over $k$ with $l=k[x]/(x^3-b)$ cubic
\'{e}tale. If the binary cubic form $\langle  a,b\rangle$
 represents $1$  or $b^2$, or if the binary cubic form $\langle  a,b^2\rangle$ represents $1$ or $b$,
 then $A$ splits, i.e. $A\cong\mbox{Mat}_3(k)$.

 More generally, take any odd integer $d$ and let $k$ be a field of characteristic not dividing $d$
  containing a primitive $d^{th}$  root of unity.
  Suppose that $A\cong (l,a)$ is a cyclic central simple algebra over $k$ of degree
 $d$ with $l=k[x]/(x^d-b)$ a field extension of $k$. Then $A$ splits, if
  the form $\langle  a,b^r\rangle$ of degree $d$ represents $b^s$
for  suitable integers $r$, $s$ with $r\not=s$, $1\leq r\leq d-1$, $0\leq s\leq d-1$.

Also for certain 27-dimensional exceptional simple cubic Jordan algebras,
 it suffices to check a much simpler 6-dimensional
cubic form  to see if the algebra is not a division algebra:
let  $A \cong (l,a)$ be a cyclic algebra over $k$ of degree 3 as above, which does not have zero divisors.
If the cubic form $\langle  a,c\rangle\otimes n_{l/k}$
 represents $1$  or $a^2$, or if the cubic form $\langle  c,a^2\rangle\otimes n_{l/k}$ represents $1$ or $a$,
 then the algebra $J=J(A,c)$ obtained from $A$ and $c\in k^\times$ via the first Tits construction is not a Jordan division algebra.

\section{Main theorem}

\subsection{} A  {\it form of degree} $d\geq 2$ over $k$
 is a homogeneous polynomial $\varphi\in k[x_1,\dots,x_m]$ of degree $d$ over $k$.
If we can write $\varphi$ in the form $\varphi(x_1,\dots,x_m)=a_1x_1^d+\ldots +a_mx_m^d$ we use the notation $\varphi=\langle  a_1,\ldots
,a_n\rangle  $ and call the form $\varphi$ {\it diagonal}. A form $\varphi$ over $k$ {\it represents} an element
$a\in k^\times$ if there are $a_i\in k$ such that $\varphi(a_1,\dots,a_m)=a$.

\subsection{}
A finite dimensional commutative $k$-algebra $l$ is called {\it \'{e}tale},
if $l\cong k_1\times\dots\times k_r$ for some finite separable field extensions $k_1,\dots, k_r$ of $k$.
A {\it Galois $G$-algebra} over $k$ is an \'{e}tale $k$-algebra $l$ endowed with an action by $G$,  where  $G$ is
 a group of $k$-automorphisms of $l$, such that the order of $G$ equals the dimension of $l$ over $k$, and such
that $l^G=\{x\in l\,|\, g(x)=x \text{ for all } g\in G \}=k$ [KMRT, (18.B.), p.~287].
A Galois $G$-algebra structure on a field $l$ exists if and only if the extension $l/k$ is Galois with Galois
group isomorphic to $G$ [KMRT, (18.16), p.~288].

 Let $l$ be a Galois $(\mathbb Z/n\mathbb Z)$-algebra over $k$. The {\it cyclic algebra} $(l,a),\, a\in
k^\times$, is defined as
$$(l,a)=l\oplus lz\oplus\ldots\oplus lz^{n-1}$$
 where $z^n=a$ and $zu=\mu (u)z$ for
$u\in l,\, \mu=1+n\mathbb Z\in \mathbb Z/n\mathbb Z$. It is a central simple $k$-algebra of degree $n$.
Every central simple algebra over $k$ of degree $n$ which contains $l$ has the form $(l,a)$ [KMRT, (30.A.), p.~414].

\begin{theorem} Let $d$ be an odd integer and let $k$ be a field such that
${\rm char}\,k$ does not divide $d$, which contains a primitive $d^{th}$ root of unity $\omega$.
Let $A=(l,a)$ be a cyclic central simple algebra over $k$ of degree $d$, such that $l=k[x]/(x^d-b)$
is a field extension of $k$ of odd degree $d$.
 If the form
$\langle  a,b^r\rangle\text{ of degree }d\text{ represents }b^s,$ for suitable integers $r$ and $s$, where
$1\leq r\leq d-1$, $0\leq s\leq d-1$ and $r\not=s$,
 then $A\cong\mbox{Mat}_d(k)$.
\ignore{(i) If the form $$\langle  a,b^{d-r}\rangle\text{ of degree }d\text{ represents }b^{d-r-1},$$
for a suitable integer $r$ with
 $1\leq r\leq d-1$, then $A\cong\mbox{Mat}_d(k)$.\\
(ii) If the form
$\langle  a,b^r\rangle \text{ of degree } d \text{ represents }1,$
 for a suitable integer $r$ with
 $1\leq r\leq d-1$, then $A\cong\mbox{Mat}_d(k)$.\\}
\end{theorem}

\begin{proof}  Let  $l=k[x]/(x^d-b)$
be a field extension of $k$ of degree $d$. Write $l=k(\alpha)$ with $\alpha^d=b$.
Let $u_1+u_2\alpha+\dots+u_{d}\alpha^{d-1}$ be an arbitrary element in $l$. Then
$$n_{l/k}(u_1+u_2\alpha+\dots+u_{d}\alpha^{d-1})=u_1^d+bu_2^d+\dots+b^{d-1}u_{d}^d+\text{ sum of ``mixed terms''}$$
for $u_i\in k$, where ``mixed terms'' are terms of the kind $u_{i_1}u_{i_2}\cdots u_{i_{d}}$ with not all $i_j$ identical.

 Let $r$ be an integer such that $1\leq r\leq d-1$.
 If the form $\langle  a,b^r\rangle$ of degree $d$ represents
$b^s$ for a suitable integer $s$ with $0\leq s\leq d-1$ and $r\not=s$,
and if $b^s$ is not represented by the form $\langle  b^r\rangle$, then there are $x,y\in k$ such that
$$ax^d+b^ry^d=b^s$$
 and we must have $x\not=0$. Thus $a=b^s(1/x)^d+b^r(-y/x)^d=n_{l/k}((1/x) \alpha^s+ (-y/x)\alpha^r)$
and 
we obtain the assumption. It remains to see that indeed $b^s$ is not represented by the form $\langle  b^r\rangle$.

Suppose $b^s=b^r\gamma^d$ for some $\gamma\in k$ is indeed possible. Then $b^s/b^r=b^{r'}=\gamma^d$ with $r'=s-r\not=0$ and
 $\gamma^d=(\alpha^{r'})^d$, hence $(\alpha^{r'}/\gamma)^d=1$.
Therefore  $\alpha^{r'}/\gamma\in k(\alpha)$ is a $d^{th}$ root of unity. By assumption, $k$ contains all primitive
$d^{th}$ roots of unity, thus $\alpha^{r'}/\gamma\in k$. This implies that $\alpha^{r'}\in k$, so that the degree of
$\alpha$  over $k$ is less than or equal to $r'$ with
${r'}\leq d-1$, a contradiction.
\ignore{
If the form $\langle  a,b^r\rangle  $ represents $1=b^0$,
there are $x,y\in k$ such that $$ax^d+b^ry^d=1.$$
If $x=0$ it follows that $b^r=\frac{1}{y^d}$.

Suppose $b^r=\gamma^d$ for some $\gamma\in k$ is indeed possible.
 Then $\gamma^d=(\alpha^r)^d$, hence $(\alpha^r/\gamma)^d=1$.
Therefore $\alpha^r/\gamma\in k(\alpha)$ is a $d^{th}$ root of unity, thus $\alpha^r/\gamma\in k$.
This implies that $\alpha^r\in k$, so that $\alpha$ has degree
$\leq r\leq d-1$ over $k$, a contradiction.
}
\end{proof}

The question whether or not there might be other binary forms of degree $d$ which yield similar sufficient conditions
is left open here. Theorem 1 becomes the following result for $d= 3$:

\begin{corollary} Let $k$ be a field of ${\rm char}\,k\not=3$ containing a primitive third root of unity $\omega$.
 Let $A=(l,a)$ be a cyclic
central simple algebra over $k$ of degree $3$, where
$l=k(\alpha)$ with $\alpha^3=b$ is a cubic field extension of $k$. If
 the binary cubic form $$\langle  a,b\rangle \text{ represents 1 or }b^2, $$ or
 the binary cubic form
$$\langle  a,b^2\rangle  \text{ represents 1 or }b,$$
then $A\cong\mbox{Mat}_3(k)$.
\end{corollary}

\begin{remark}\ignore{
 (i) If one also checks diagonal forms of degree $d$ and dimension greater than two (i.e. forms in more than two variables),
one gets more sufficient conditions for the algebra to be split; e.g,
let $d=3$ and let $1$ be represented by the cubic form $\langle a,b,b^2\rangle$, but not
by the binary cubic form $\langle b,b^2\rangle$. Then $$1=ax^3+
by^3+b^2z^3$$ and we must have $x\not=0$. Thus $a=\frac{1}{x^3}(-by^3-b^2 z^3)=b(-\frac{y}{x})^3+b^2(-\frac{z}{x})^3
=n_{l/k}(-\frac{y}{x}\alpha-\frac{z}{x}\alpha^{2})$ and again $A\cong\mbox{Mat}_d(k)$.
}
(i) Let $k=\mathbb{Q}_p$ be the $p$-adic numbers, $p\not=3$, and let $a=1$,
$b=p$. Then $1$ is represented by the binary cubic form $\langle 1,p\rangle$, but $p^2$ is not.
\ignore{ if $b^2$ is represented by
$\langle a,b\rangle$, then $\langle a,b, -b^2\rangle= \langle 1,p,-p^2\rangle$ must be isotropic. However,
 by Springer's theorem [Mo] it is anisotropic, since the three residue forms $\langle 1\rangle$,
 $\langle 1\rangle$ and $\langle -1\rangle=\langle 1\rangle$ are all anisotropic.}
  Thus the first two conditions in Corollary 1 are not equivalent.\\
(ii) Corollary 1 holds even when $l$ is not a field extension: if you assume that $l=k[x]/(x^3-b)$ and that
$l$ is not a field extension, then $b=\alpha^3$ with $\alpha\in k$, and
$l\cong k\times k \times k$. Thus, the algebra $(l,a)$ is split for every $a\in k^\times$.
Corollary 1 becomes a tautology in this case.
\ignore{
(iii) By a well-known theorem of Wedderburn, every central simple algebra of degree 3 over $k$ is cyclic.
 By [KMRT, (18.32), p.~300],
every cubic  \'{e}tale algebra $l$ over $k$, ${\rm char}\,k\not=3$, is isomorphic to an algebra of the kind
$k[x]/(f)$, where we may choose the polynomial $f$ to be of the above form $f=x^3-b$ for some $b\in k^\times$, if and only
if the discriminant algebra of $l$ has the form $\Delta(l)\cong k[t]/(t^2+t+1)$.
}
\end{remark}

\begin{example} (communicated by P. Morandi.)
There are cyclic algebras which are split over $k$, but do not satisfy any of the three
conditions given in Corollary 1.

Let $k$ be a field of ${\rm char}\,k\not=3$ containing a primitive third root of unity $\omega$ and let
$F=k(x,y,z,t)$ be the rational function field in four variables over $k$. We
define $a=x^{3}+y^{3}t+z^{3}t^{2}-3xyzt$ and $b=t$.
Consider the cyclic algebra
$$A=\left(  F\left(\sqrt[3]{t}\right)  ,a\right), $$
 then $A$ is split since
$$n_{F\left(\sqrt[3]{t}\right)  /F}(x+y\sqrt[3]{t}+z\left(  \sqrt[3]{t}\right)  ^{2})=a.$$
However, none of the conditions of Proposition 1 hold. To see this, first
suppose that $\left\langle a,t\right\rangle $ represents $1$. Then there are
$u,v,w\in k[x,y,z,t]$ with $w\neq0$ and $au^{3}+tv^{3}=w^{3}$. This is
impossible considering degrees in $t$ modulo $3$; we have $\deg_{t}%
(au^{3})\equiv2\operatorname{mod}3$, $\deg_{t}(tv^{3})\equiv
1\operatorname{mod}3$, and $\deg_{t}(w^{3})\equiv0\operatorname{mod}3$. Next,
suppose that $\left\langle a,t\right\rangle $ represents $t^{2}$. Then there
are $u,v,w\in k[x,y,z,t]$ with $\gcd(u,v,w)=1$ and $au^{3}+tv^{3}=t^{2}w^{3}$.
This equation implies $t$ divides $u$ since $t$ does not divide $a$. If
$u=tu_{0}$ for some polynomial $a_{0}$, then we obtain $at^{2}u_{0}^{3}%
+v^{3}=tw^{3}$, which implies $t$ divides $v$. Writing $v=tv_{0}$ and
substituting, we see that $t$ divides $w$. This contradicts $\gcd(u,v,w)=1$.
Finally, suppose that $\left\langle a,t^{2}\right\rangle $ represents $1$.
Then there are $u,v,w\in k[x,y,z,t]$ with $w\neq0$ and $au^{3}+t^{2}%
v^{3}=w^{3}$. If we specialize $x=0$, we obtain an equation of the form
$(ty^{3}+t^{2}z^{3})p^{3}+t^{2}q^{3}=r^{3}$ for some $p,q,r\in k[y,z,t]$ with
$\gcd(p,q,r)=1$. We argue as in the previous case. The equation shows that $t$
divides $r$. Writing $r=tr_{0}$, we get $(y^{3}+tz^{3})p^{3}+tq^{3}=t^{2}%
r_{0}^{3}$, which implies that $t$ divides $p$. Setting $p=tp_{0}$ and
substituting, we get $(y^{3}+tz^{3})t^{2}p_{0}^{3}+q^{3}=t^{2}r_{0}^{3}$,
which implies that $t$ divides $q$, contradicting the assumption
$\gcd(p,q,r)=1$.
Thus, the converse to Theorem 1 is false.
\end{example}

\subsection{}
The above idea to restrict the norm of an algebra to suitably chosen subspaces of the algebra  can also be applied to certain nonassociative algebras to obtain sufficient
conditions for them to be not a division algebra. Let us recall a way to construct 27-dimensional
exceptional simple cubic Jordan algebras (also called {\it Albert algebras})
out of central simple  associative algebras of degree 3 first: Let $A$ be a  central simple associative
algebra over $k$ of degree 3 with reduced norm $n_{A/k}$,
reduced trace $tr_{A/k}$ and adjoint $\sharp$. The vector space $V=A\oplus A\oplus A$ together with
$$\begin{array}{l}
 N(a_1,a_2,a_3)=n_{A/k}(a_1)+cn_{A/k}(a_2)+c^{-1}n_{A/k}(a_3)-tr_{A/k}(a_1a_2a_3),\\
(a_1,a_2,a_3)^\sharp=(a_1^\sharp-a_2a_3,c^{-1}a_3^\sharp-a_1a_2,ca_2^\sharp a_3a_1),\\
T((a_1,a_2,a_3),(b_1,b_2,b_3))=tr_{A/k}(a_1b_1)+tr_{A/k}(a_2b_3)+tr_{A/k}(a_3b_2).
\end{array}$$
becomes a central simple exceptional Jordan algebra with identity $(1,0,0)$ and $U$-operator
$U_vw=tr_{A/k}(v,w)v-v^\sharp \times w$.
(As customary, we write $((a_1,a_2,a_3),(b_1,b_2,b_3))$ $\to (a_1,a_2,a_3)\times (b_1,b_2,b_3)$ for the bilinearization
of the quadratic map $\sharp$.)  It is denoted by $J(A,c)$ and called a {\it first Tits construction} (for a more detailed explanation
on the terminology and the construction, see, for instance,
 [KMRT, p.~525]).

\begin{theorem} Let $k$ be a field of ${\rm char}\,k\not=3$ containing a primitive third root of unity $\omega$.
 Let $A=(l,a)$ be a cyclic
central simple division algebra over $k$ of degree $3$, where
$l=k(\alpha)$ with $\alpha^3=b$ is a cubic field extension of $k$ and $c\in k^\times$. If
 the cubic form $$\langle  c,a\rangle\otimes n_{l/k} \text{ represents 1 or }a^2, $$ or
 the cubic form
$$\langle  c,a^2\rangle \otimes n_{l/k} \text{ represents 1 or }a,$$
then the Albert algebra $J=J(A,c)$ is not a Jordan division algebra.
\end{theorem}

\begin{proof} The exceptional simple Jordan algebra $J=J(A,c)$ is a division
 algebra if and only if $c$ is not a norm of $A$ [M, Theorem 6].
Assume that $\langle  c,a\rangle\otimes n_{l/k}$ represents 1.
Therefore there exist $x,y\in l$ such that $1=cn_{l/k}(x)+a n_{l/k}(y)$. If $x=0$, then $1=a n_{l/k}(y)$ implies
that $a$ is a norm of $l$, contradicting our assumption that $A$ is a Jordan division algebra. Hence we conclude from
$c= n_{l/k}(\frac{1}{x})+a n_{l/k}(\frac{-y}{x})$ that $c$ is a norm of $A$.
The rest of the assertion is shown analogously.
\end{proof}

\smallskip
{\it Acknowledgements:}  The author would like to thank P. Morandi for providing a counterexample and
 would like to acknowledge financial support of the ``Georg-Thieme-Ged\"{a}chtnisstiftung''
 (Deutsche Forschungsgemeinschaft) during her stay at the University of Trento.
 She thanks the Department of Mathematics at Trento for its hospitality and congenial atmosphere.

\end{document}